\documentclass[a4paper,12pt]{article}
\usepackage{ucs}
\usepackage{amsfonts} 
\usepackage{amsmath} 
\usepackage{graphicx,subfigure}
\usepackage{caption} 
\usepackage{sidecap} 
\usepackage{enumerate}
\usepackage{MnSymbol}
\usepackage{color}
\usepackage{xfrac}

\title{}
\author{}
\date{}

\newtheorem{theorem}{Theorem}[section]
\newtheorem{lemma}[theorem]{Lemma}
\newtheorem{proposition}[theorem]{Proposition}

\newenvironment{keywords}{\begin{@abssec}{\keywordsname}}{\end{@abssec}}
\newenvironment{@abssec}[1]{%
\if@twocolumn
\section*{#1}%
\else
\vspace{.05in}\footnotesize
\parindent .2in
{\upshape\bfseries #1. }\ignorespaces
\fi}
{\if@twocolumn\else\par\vspace{.1in}\fi}
\newenvironment{AMS}{\begin{@abssec}{\AMSname}}{\end{@abssec}}

\newcommand{\qed}{\nobreak \ifvmode \relax \else
\ifdim\lastskip<1.5em \hskip-\lastskip
\hskip1.5em plus0em minus0.5em \fi \nobreak
\vrule height0.75em width0.5em depth0.25em\fi}

\newcommand\keywordsname{Key words}
\newcommand\AMSname{AMS subject classifications}

\newcommand{\allgreen}{\color{green}{}}

\newcommand{\Saddle}{\mathcal{P}}
\newcommand{\Repel}{\mathcal{R}}
\newcommand{\Cone}{\mathcal{C}} 
 
\newcommand{\upPhi}{\hat{\Phi}}
\newcommand{\lowPhi}{\Phi}
\newcommand{\Torus}{\mathbb{T}^2} 
\newcommand{\torus}{\mathbb{T}} 
\newcommand{\TorusD}{\mathbb{T}^d} 
 
\newcommand{\TorusDmO}{\mathbb{T}^{d-1}}
\newcommand{\Tile}{P}
\newcommand{\Plane}{\mathcal{P}}
\newcommand{\Man}{\mathcal{M}}
\newcommand{\HDC}{UDV\ }
\newcommand{\HDCp}{UDV}
\newcommand{\R}{\mathbb{R}}

\usepackage{bm}

\begin{document}
\title{Multi-chaos from Quasiperiodicity}
\author{Suddhasattwa Das\footnotemark[1], \and James A Yorke\footnotemark[3]}
\footnotetext[1]{Courant Institute of Mathematical Sciences, New York University}
\footnotetext[3]{University of Maryland, College Park}
\date{\today}
\maketitle

\begin{abstract}
One of the common characteristics of chaotic maps or flows in high dimensions is ``unstable dimensional variability'', in which there are periodic points whose unstable manifolds have different dimensions. In this paper, in trying to characterize such systems we define a property called ``multi-chaos''. 
A set $X$ is multi-chaotic if $X$ has a dense trajectory and for at least 2 values of $k$, the $k$-dimensionally unstable periodic points are dense in $X$. 
All proofs that such a behavior holds have been based on 
hyperbolicity in the sense that (i) there is a chaotic set $X$ with a dense trajectory and (ii) in X there are two or more hyperbolic sets with different unstable dimensions.
We present a simple 2-dimensional paradigm for multi-chaos in which a quasiperiodic orbit plays the key role, replacing the large hyperbolic set.
\end{abstract}



\section{Introduction} \label{sect:intro}
The type of chaos seen in the most common examples is low-dimensional and not representative of typical systems with several positive Lyapunov exponents. We introduce a term ``multi-chaos'' to describe certain chaotic sets in order to emphasize how different chaotic systems can be from the low dimensional systems.

{\bf A definition of chaos.} There are a variety of definitions of chaos (see Sander and Yorke \cite{sander:yorke:15} and Hunt and Ott \cite{ChaosDef}) 
and we have chosen a definition below that fits the needs of this project. Let $M$ be a smooth manifold and
$F:M\to M$ is a smooth map. We say  a set $X\subset M$ is {\bf transitive} if
(i) $X$ is an uncountable closed set;
(ii) $X$ is invariant ($F(X)=X$); and 
(iii) $X$ has a dense trajectory $F^n(x_0)$ where $n=1,2,\cdots$ for some  $x_0\in X$. 
We will say such an $x_0$ is a {\bf transitive point} (of $X$). 

We say $X$ is {\bf chaotic} if $X$ is a transitive set and periodic points of $F$ are dense in $X$.

{\bf Unstable Dimension Variability (\HDCp).} We say a chaotic set $X$ has 
\emph{unstable Dimension Variability} or is \emph{Unstable-Dimension Variable (\HDCp)} if there are two periodic orbits in $X$ whose unstable manifolds have different dimensions. This term was introduced in \cite{kostelich97} without a precise definition and has been used frequently since then.

Such behavior is quite different from the most common examples of chaos, but the first examples are quite old.
Abraham and Smale \cite{Omega_stab} introduced a diffeomorphism in four dimensions 
that they showed was robustly non-hyperbolic and Simon \cite{simon72} reduced the required dimension to three. 
Bonatti and Diaz \cite{Blender2, Blender3} richly expanded the family of robust examples.
Here ``robust'' means that all sufficiently small perturbations of the system preserve the structure. 
That may be rephrased as saying the set of such dynamical systems is an open set in the space of all dynamical systems. We sometimes call this concept of robustness ``topological'' robustness to distinguish from other properties that might be ``probabilistically'' or measure-theoretically robust.
Their examples had periodic orbits of different unstable dimensions in a single chaotic set so they was \HDCp.
The first robust dynamical systems were those that had hyperbolic dynamics, and these examples greatly expanded the types of robust behaviors.
%

{\bf Multi-chaos.}  Let $X$ be chaotic.
We say $X$ is $k$\textbf{\textit{-chaotic}} for some integer $k>0$ if there is a dense set of periodic points in $X$ whose unstable manifolds intersected with $X$ have (topological) dimension $k$. More loosely we can say there is 
$k$\textbf{\textit{-chaos.}} 

We will say $X$ is \textbf{\textit{multi-chaotic}} if 
$X$ is $k$-chaotic for two or more values of $k$. We sometimes write the map $F$ is multi-chaotic on (the set) $X$, or we can say the map has 
\textbf{\textit{multi-chaos}}. For example if $X=M$ is a two dimensional torus and is chaotic, then it is multi-chaotic if it is $1$-chaotic and $2$-chaotic. That is, it has a dense set of periodic points that are saddles and another dense set of periodic points that are repellers. This is the type of situation we will construct in detail. 
We assume that the definitions of ``saddle'' and ``repeller''  require the orbits to be hyperbolic.
The goal of this paper is to describe new simple, verifiable conditions which lead to multi-chaos. Our main theorem is stated in dimension 2, but it provides a new paradigm for multi-chaos in any dimension $\geq 2$. 

{\bf A \HDCp:multi-chaos conjecture.} Notice that a \HDC map must have at least two periodic orbits while a multi-chaos map must have two corresponding dense sets of periodic orbits. 
So it appears that multi-chaos is much more restrictive. None the less we know of no \HDC that is not multi-chaotic.
We conjecture that almost every (in the sense of prevalence \cite{hsy92}) \HDC map is multi-chaotic. 
We note that in the sense of prevalence \cite{hsy92}, almost every smooth map has the property that all of its periodic orbits are hyperbolic \cite{hsy92}.

\textbf{Quasiperiodicity.} 
Our main result gives examples of multi-chaotic maps for which there is a quasiperiodic trajectory.
Recall two dynamical systems are said to be \textbf{conjugate} if there is a continuous change of coordinates transforming one to the other in which the change of coordinates has a continuous inverse. 
This change of variables is called a \textbf{conjugacy}. 
A map $h$ on a circle $S^1$ or torus $\TorusD$ which is conjugate to a rotation $\theta\mapsto\theta+\rho\bmod 1$ where $\rho$ is irrational is said to be {\bf quasiperiodic}.

{\bf An example of multi-chaos from quasiperiodicity.} To show how multi-chaos arises from quasiperiodicity, we examine a rather simple version of our main result, Theorem \ref{thm:Main}.
Expressions like $(x,y) \bmod 1$ imply $\mod1$ is applied in each coordinate. Consider

\begin{equation}\label{eqn:map}
F_0(x,y) = (mx,\ ax + y+g_0(x,y))\ \bmod 1,
\end{equation}
where $a$ and $m$ are integers with $|m|>1$, and $g_0:\mathbb{R}^2\rightarrow\mathbb{R}$ is $C^1$ and $\mathit{\mathbb{Z}}$\textbf{\textit{-periodic}} in $x$ and $y$; that is, $g_0(x+n_1,y+n_2) = g_0(x,y)$ for all integers $n_1$ and $n_2$. 
These maps are called \textbf{skew-product} maps. Our main result, Theorem \ref{thm:Main}, however is for a more general form of torus map that  contain a quasiperiodic curve, from which we conclude there is multi-chaos.  Also
$G$ can be much larger.
%
The map has the following dynamics in the X coordinate.
\begin{equation}\label{eqn:expnd_circ}
x\mapsto mx\bmod 1.
\end{equation}

A \textbf{vertical circle} is the set $S_{x_0}:=\{(x,y)\in\Torus:x=x_0\}$. The map (\ref{eqn:map})  above takes vertical circles onto vertical circles, that is, for every $x\in S^1$,
\begin{equation}\label{eqn:vertical}
F_0(S_x)=S_{x^*},\ \mbox{ where }x^*:=mx\bmod 1
\end{equation}
For each $x_0$ that is periodic under the map in Eq. \ref{eqn:expnd_circ} with period $n$, we have $F_0^n(S_{x_0})=S_{x_0}$. These will be called the \textbf{periodic circles} of the map. Each may or may not be quasiperiodic. Notice that if the map $F_0^n$ has a quasiperiodic curve, then it must be a vertical circle. It must have a constant value of x.

Let $F:\TorusD\to\TorusD$ be continuous, where for each 
$x\in\R^d$, we can write $x\bmod1\in \torus$, applying mod 1 to each coordinate.
Then $F$ can be uniquely decomposed into a linear part and a periodic part.
That is, there is a $d\times d$ matrix $M$ with integer entries 
and $G:\mathbb{R}^d\to\mathbb{R}^d$ is a continuous periodic function with period 1 in each coordinate, such that (abusing terminology a bit)
\begin{equation}\label{eqn:periodic_part}
F(z)=[Mz+G(z)]\bmod1\mbox{ for all }z\in\Torus.
\end{equation}
Note that if $\nu$ is a d-dimensional vector with integer entries, then
$G(z+\nu)=G(z)$ while $A(z+\nu)-Az$ is a vector with integer entries so
$F(z+\nu) = F(z)$ for all $z$ due to the ``mod 1'' in Eq. \ref{eqn:periodic_part}.

Our example uses the hypotheses: \\

\noindent
(${\bm E_{M}}$): The integer matrix $M$ has  two integer eigenvalues, $1$ and $m$ with $|m| > 1$.
\\(${\bm E_{\Saddle\Repel}}$): $F$ has a periodic saddle $\Saddle$ and a periodic repeller $\Repel$.
\\(${\bm E_{QuasiP}}$): There is a curve on which $F^p$  is quasiperiodic for some integer $p$.\\

\noindent
{\bf Example A.} \emph{
Let $F:\Torus\to\Torus$ be a $C^1$ map of the form
  $F(z)=Mz+G(z) \bmod1$  where the matrix $M$ satisfies ($E_{M}$) and $G$ is periodic with period $1$ in each coordinate. Assume ($E_{\Saddle\Repel}$) and ($E_{QuasiP}$).
Then there is a $\delta(M) >0$ such that if the function $G$  satisfies $\|DG\|<\delta(M)$, then the torus $\Torus$ is multi-chaotic for $F$. }
\bigskip

We state an analogue of the above example involving diffeomorphisms, and leave it to the reader to adapt the proofs to this case. We use the following condition.

\bigskip
\noindent
(${\bm E_{M}^3}$): The $3\times3$ integer matrix $M$ has  three eigenvalues, $1$ and $m$ and $1/m$ with $|m| > 1$.\\

\noindent
{\bf Example B.} \emph{
Let $F:\torus^3\to\torus^3$ be a $C^1$ diffeomorphism of the form
  $F(z)=Mz+G(z) \bmod1$  where the $3\times3$ matrix $M$ satisfies ($E_{M}^3$) and $G:\mathbb{R}^3\to\mathbb{R}^3$ is periodic with period $1$ in each coordinate. Assume ($E_{\Saddle\Repel}$) and ($E_{QuasiP}$).
Then there is a $\delta(M) >0$ such that if the function $G$  satisfies $\|DG\|<\delta(M)$, then the torus $\torus^3$ is multi-chaotic for $F$. }
\bigskip

\textbf{Countably infinitely many invariant curves.} 
Example A and Theorem \ref{thm:Main} rely on the existence of a quasiperiodic curve, a curve which is invariant under $F^n$ for some $n$.
As the proof reveals,  
there is a one-to-one correspondence between the period-$n$ curves of $F$ in Example A and the period-$n$ circles of Eq. \ref{eqn:map}.

Specifically, by Eq. \ref{eqn:vertical}, for every $k,n\in\mathbb{N}$, the vertical circles $S_x$ with $x=k/(m^n-1)$ corresponds to an invariant curve under $F^n$. 
Assume in this paragraph $F$ is $C^\infty$.
Assume the torus map is smoothly and generically
parameterized by some parameter $t\in\R$. Each periodic curve varies continuously with $t$. 
For each t, each periodic curve (corresponding to some $k/(m^n-1)$) will either have a periodic orbit on it or will be quasiperiodic. The set of $t$ for which the curve has a periodic attracting orbit (and hence is not quasiperiodic) is typically an open dense set. 

The intersection $J$ of this countable collection of open dense sets is dense in $\R$.
The complement of $J$ is the set of t for which there is at least one quasiperiodic curve. Our numerical calculations with Yoshi Saiki seem to suggest that this is a full measure set in some interval of $t$ values; that is for some $t_1$ and $t_2>t_1$, the quasiperiodic set of parameters  has full measure in $[t_1,t_2]$.

{\bf Topological versus measure-theoretic robustness.} 
One example is instructive.
Smale and Williams \cite{smale-williams-76} analyzed the logistic map $f_b(x) = bx(1-x)$ for the case $b = 3.83$, where there is a period-three attracting periodic orbit. They showed that $f_b$ is robust. They use the term ``structurally stable'' instead of ``robust''. Each $b$ has at most one attractor, and $b = 3.83$ also has a chaotic hyperbolic set lying on boundary of the basin of attraction, and these invariant sets change continuously in the interval where there is a period-three attractor, that is, in the period-3 ``window''. If we examined only the $b$ for which $f_b$ is robust, we would miss all of the $b$ for which there are chaotic attractors. Studying only robust dynamical systems seems to have its negative aspects. 

There is a chaotic invariant set for each $b\in C_0 := (3.57,4]$, including points where there is a periodic attractor and at least one chaotic saddle.
There are infinitely many periodic windows where there is a periodic attractor and these are dense in $C_0$ (Milnor and Thurston \cite{milnor-thurston-umd}). Jacobson \cite{jakobson80} showed that the complement $C\subset C_0$ of the periodic windows in the chaotic interval $C_0$ is a set of positive measure, so the sum of the lengths of these topologically robust intervals is smaller than the length of $C_0$. 
That set could be called ``measure-theoretically'' robust but the topologically robust intervals taken together are small. 

Since C has positive Lebesgue measure, almost-every $b\in C$ has the property that in each sufficiently small neighborhood of $b$ only a small fraction of the interval will be in a periodic window. 
So while such periodic intervals are dense, locally they are a small fraction.

\begin{figure}
\includegraphics[width=.45\textwidth]{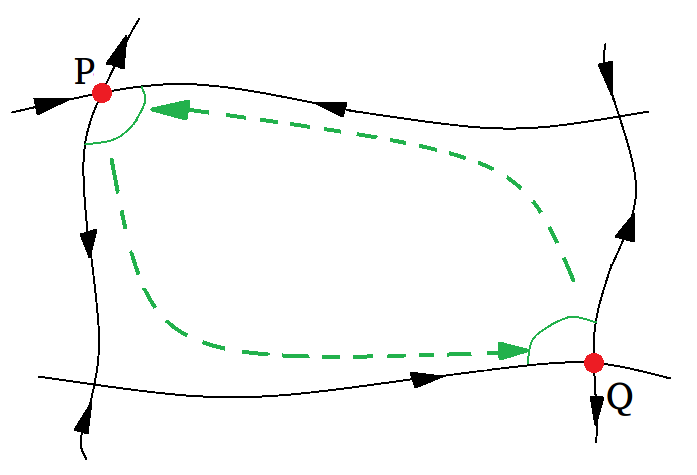}
\includegraphics[width=.45\textwidth]{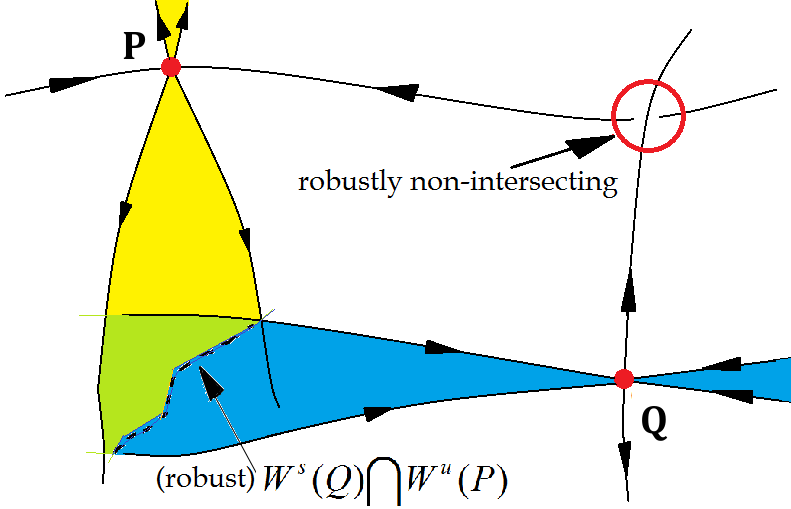}
\label{fig:Variable_unstable_dimension_problem}
\caption[Heteroclinic cycles]{A Heteroclinic cycle is created when stable and unstable manifolds of a saddle $P$ transversely intersect the unstable and stable manifolds respectively of another saddle $Q$. Part(a) displays such a case for two saddles in a 2-dimensional setting. Since the unstable dimensions of $P$ and $Q$ are the same, the two intersections persist under small perturbations of the dynamical system. 
The regions marked in green are portions of neighborhoods of $P$ and $Q$ and each of these regions have a Cantor set of points which keep coming close to $P$ and $Q$ infinitely many times in their trajectories. Part(b) shows a configuration in a 3-dimensional setting in which the unstable dimensions of $P$ and $Q$ are $2$ and $1$ respectively. 
As a result, for almost every perturbation of the system, $W^s(P)$ and $W^u(Q)$ will not intersect. 
Thus, the configuration of a heteroclinic cycle is absent and some other mechanism is needed by which points near $P$ can land close to $Q$ and vice versa. 
}
\label{fig:Unstable_dimension_illustration}
\end{figure}

\textbf{The \HDC problem.}
Consider a dynamical system with two hyperbolic periodic points $P$ and $Q$, as shown in  Fig. \ref{fig:Unstable_dimension_illustration}. 
If they have the same unstable dimension and in addition, the stable manifold of each intersects the unstable manifold of the other, the configuration is a heteroclinic cycle. 
Then there is an invariant Cantor set with a dense trajectory that contains both $P$ and $Q$. 
Hence the dense trajectory comes arbitrarily close to both $P$ and $Q$. If the unstable dimension of $Q$ is more than that of $P$, then generically, the stable and unstable manifolds of $P$ and $Q$ respectively will not intersect since the sum of their dimensions is less than the dimension of the ambient space. Therefore, in such a case, the mechanism of heteroclinic intersections is lost. 
So we are faced with the problem of determining why should there be trajectories arbitrarily near $P$ which pass arbitrarily close to $Q$ and return close to $P$, etc., a necessary feature of chaotic sets with varying unstable dimensions.

\textbf{A consequences of \HDCp: Non-Shadowability.} Yuan and Yorke \cite{NonShadow}, C. Grebogi et. al. \cite{NonShadow2} proved the property of  ``non-shadowability'' for each family of maps that has a recurrent set containing two saddles of different unstable dimensions. Hence the chaotic set of a \HDC map never has the shadowing property. Several other papers, like Abdenur and Diaz \cite{diaz2007},  discuss the consequences of co-existence of saddles of different unstable dimensions in the same chaotic set.

\textbf{A consequence of \HDCp: Oscillation of finite time Lyapunov exponents (FTLE).} The transitive trajectory comes arbitrarily close to each of the two periodic sets and therefore on rare occasions it spends arbitrarily long intervals of time near each of the orbits. 
Suppose FTLE of a transitive trajectory are computed using the trajectory's iterates $n+1$ to $n+N$. Let  $\mathcal N_{n+1}^{n+N}$ denote the number of positive FTLEs for that segment of the trajectory from $n+1$ to $n+N$. 
When the trajectory is extremely close to a periodic orbit, the FTLEs will be quite close to the Lyapunov exponents of the periodic orbit. That means that for each $N$ the number $\mathcal N_{n+1}^{n+N}$ of positive finite-time Lyapunov exponents will fluctuate. In other words, {\bf one of the finite time Lyapunov exponents oscillates about zero}.
See Dawson et al \cite{dawson94}. This can be thought of as a test for \HDCp.

{\bf Examples of \HDCp.} 
Smale and Abraham\cite{Omega_stab} constructed the first robust \HDC example. Bonatti and Diaz \cite{Blender2, Blender3} introduced the concept of \textbf{\textit{blenders}} to generalize a hyperbolic set that occurs in the earlier \HDC examples enabling them to construct additional examples of robust diffeomorphisms is \HDCp. 
By the definitions of Bonatti and Diaz, the earlier examples had blenders but the definition of blenders also allows new examples.
Also see \cite{Blender1, Blender4} for other constructions using blenders.

Diaz and Gorodetski \cite{UnstDimVar1} constructed a residual set of diffeomorphisms such that each homoclinic class that contains saddles of different ``indices'' (the number of unstable dimensions) also contains an uncountable support of an invariant ergodic non-hyperbolic measure of $f$. By smoothing piecewise continuous skew-product maps \cite{Horse_solen}, Gorodetski and Ilyashenko \cite{Robust_inv, Regular_centre} constructed \HDC  diffeomorphisms in a locally maximal invariant set. 

 \textbf{\HDC in maps with blow-out bifurcations.} In numerical studies of ``blow-out bifurcations'' and ``riddled basins'', one often sees chaotic sets that contain a piece of an invariant plane of dimension $k_1$ in which there are periodic orbits of different unstable dimensions, some of which are larger than $k_1$. 
 Hence, at least part of their unstable manifolds are transverse to the plane. Often there is a larger chaotic set that contains that piece of invariant plane. Such sets have the potential of being multi-chaotic, but no examples have been rigorously shown to have that property.

\textbf{\HDC in numerical and physical experiments.} We believe that multi-chaos is present in most high-dimensional chaotic sets. Experiments on higher dimensional chaotic sets have strongly suggested that the periodic orbits frequently have different numbers of unstable dimensions. ``\HDCp''  is a closely related topic (usually studied numerically) in which the number of local positive exponents fluctuates along a typical trajectory, ``local'' meaning that the exponents are approximated using a fixed segment length of a trajectory. 

S Ponce Dawson \cite{dawson96} describes a situation analogous to blenders except that instead of using a topologically robust hyperbolic set, she uses a H\'enon attractor, which plausibly might be measure-theoretically robust. See also \cite{dawson97}. 
These papers give a valuable but non rigorous alternative view of \HDCp.

Pikovsky and Grassberger \cite{PikoGrass} and later Glendinning \cite{Milnor_attractor} examined a system of coupled identical dynamical systems and proved the existence of a two-dimensional chaotic attractor $X$ and that repeller periodic orbits are dense in the chaotic attractor $X$. Our numerical investigations with Y. Saiki strongly suggest that the system is multi-chaotic on $X$. E. Barreto and P. So \cite{Barreto} simulated coupled dissimilar chaotic attractors and found variability in the number of ``finite-time'' positive (unstable) Lyapunov exponents (along a fixed-length of the trajectory) fluctuates as time evolves.

In a hyperbolic set all periodic points have the same unstable dimension.
All of the \HDC examples above were based on the existence and interaction of hyperbolic sets whose unstable dimensions were different. Our examples here are instead based on quasiperiodic orbits. Some of these examples might have co-existing blenders, but so far there exists no test for the existence of blenders. 

\section{A theorem obtaining multi-chaos from quasiperiodicity}\label{sec:mainTheorems}

\bigskip
We first describe the concept of a dominant expanding direction in terms of expanding cones.

\textbf{Invariant expanding cones.} The property that a map has a ``a cone system'' can be defined far more generally on any manifold, but here we will define it in a restricted context on $\TorusD$. Let $w$ be a  vector on $\TorusD$ with $|w|=1$ and $W$ a $(d-1)$-dimensional subspace such that $W$ and $w$ together span the tangent space (which is $\mathbb{R}^d$).  For a tangent vector $v$ at a point $z\in\Torus$, we will write $v=(a,b)\in\mathbb{R}\times\mathbb{R}^{d-1}$ to indicate that  $v=aw+b$ where $b\in W$. 
Let $v'=DF(z)v$ where $v'=(a',b')$ in the above notation.
Then the map is said to have an (invariant), (expanding) \textbf{cone system with respect to $w$ and $W$} if there are constants $K>1$ and $\alpha>0$ such that the following are satisfied for every point $z$ and vectors $v=(a,b)$ for which $|b|\leq\alpha\|a\|$.
\\(i) $|b'|<\alpha|a'|$, and
\\(ii) $|a'|>  K|a|$.
\\(iii) for each vector $u$ transverse to the cone, $\|DF(z)u\|<K\|u\|$.
\\At every point $z\in\TorusD$, the \boldmath $(K,\alpha)$-\textbf{cone}\unboldmath, denoted as $\Cone_{K,\alpha}(z)$ is the set of vectors with representation $(a,b)$ in the tangent space at $z$ such that $|b|<\alpha|a|$. So $F$ has a cone system in the direction $w$ if there are constants $K>1$ and $\alpha>0$ such that $\Cone_{K,\alpha}(F(z))\subset DF(z)(\Cone_{K,\alpha}(z))$ and if $(a',b')=DF(z)(a,b)$, then $|a'|>K|a|$. Since $K,\alpha$ will remain fixed, we will drop them from the notation, writing $\Cone(z)$ for $\Cone_{K,\alpha}(z)$.

\begin{SCfigure}
\includegraphics[width=0.5\textwidth]{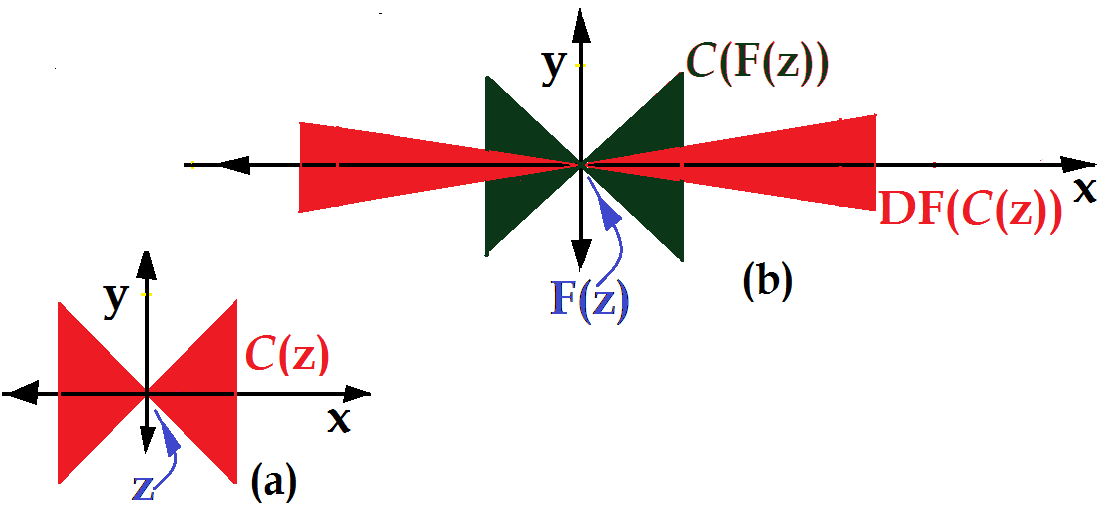}
\caption{\textbf{Cone systems.} Let $w$ be the unit vector tangent to the X axis.
Part(a) illustrates in red, a cone system in the direction $w$ at a point $z\in\Torus$, for a map $F:\Torus\to\Torus$, with $\alpha=1$. Part(b) shows in green the $\alpha$-cone at $F(z)$. Pictorially, both the $\alpha$-cones represent the vectors $\{(u,v)\ :\ |v|\leq |u|\leq 1\}$ in their respective spaces. Every vector within the red cone $\Cone(z)$ at $z$ (Part a) is mapped into the green cone $\Cone(F(z))$ at $F(z)$ (Part b) under the action of $DF(z)$; and is also stretched by a factor of at least $K>1$. Note that coordinates can be chosen so that the eigenvectors are mutually perpendicular, as is the case in this picture.}
\label{fig:inv_cone}
\end{SCfigure}

Example A illustrates the following more general theorem, which permits estimates of $\delta$, the upper bound on $\|DG\|$; note that Example A only asserts the existence of an upper bound. We introduce the following more general hypothesis:
\bigskip

\noindent
{\bf(E$_{cone}$)} There is  a $(K,\alpha)$-cone system with respect to $v_m^R, (v_m^L)^\bot$,  where $v_m^R, v_m^L$ are the right and left eigenvectors of $M$ with eigenvalue $m$.

\begin{theorem}[Main theorem]\label{thm:Main}
Let $F:\mathbb{T}^2\rightarrow\mathbb{T}^2$ be a $C^1$ map of the torus such that $DF$ is invertible everywhere. 
Assume ($E_{\Saddle\Repel}$), (E$_{cone}$), ($E_{M}$), and ($E_{QuasiP}$).
\\ Then $\Torus$ is multi-chaotic for $F$ .
\end{theorem}

\textbf{Proof of Example A as a corollary of  Theorem \ref{thm:Main}.} Example A  assumes ($E_{M}$), ($E_{\Saddle\Repel}$) and  ($E_{QuasiP}$) of  Theorem \ref{thm:Main}. 
We can choose coordinates such that $M$ takes the form $\left(\begin{array}{cc}1 & 0\\ a & m\\\end{array}\right)$, where $a$ and $m$ are integers. Then the left and right eigenvectors of $M$ corresponding to eigenvalue $m$ are  $v_m^L=(0,1)$ and $v_m^R=(am,m-1)^T$ respectively.  
Therefore, the linear map $M$ of Example A on $\Torus$ has a cone system with respect to $v_m^R$ and $(v_m^L)^\bot$, with $K=m$ and $\alpha$ any value $<\tan(\theta)$, where $\theta$ is the angle between $v_m^R$ and $v_1^R$, the right-eigenvector of $M$ with eigenvalue $1$. Having a cone system is a $C^1$ open condition, so for a periodic perturbation $g_0$ such that $\|g_0\|_{C^1}$ is small, the cone system still persists and contains the vector field $v_m$. 
Therefore, assumption (E$_{cone}$) of  Theorem \ref{thm:Main} is also satisfied for $\delta$ sufficiently small, and since $M$ is an invertible matrix, so is $F$.
Therefore the conclusion of   Theorem \ref{thm:Main} holds, i.e., $\Torus$ is multi-chaotic. 

Notice that instead of Example A's assumption that the nonlinear term is sufficiently small, the cone condition gives a verifiable condition. So for example when $a=0$, it is sufficient to take $\delta=0.5|m|$ in the modified coordinate system.

Theorem \ref{thm:Main} is proved in  Section \ref{sect:proof-Thm1}. 
The proof requires the following Conjugacy Theorem on the torus which is proved later in  Section \ref{sect:proof-conjug}.

\begin{proposition}[Conjugacy result]\label{thm:conj-exp-cone}
Let $F:\TorusD\to\TorusD$ be a $C^1$ map of the torus such that $DF$ is invertible everywhere and  $F$ satisfies assumptions ($E_{M}$), (E$_{cone}$). Then $F$ is conjugate to a continuous map of the form 
\begin{equation}\label{eqn:skew_torus}
F_0(x,y)=\left(mx, g(x,y)\right)\bmod 1,\ \mbox{ for }x\in S^1, y\in \mathbb{T}^{d-1}.
\end{equation}
i.e., there is a continuous, invertible map $H:\TorusD\to\TorusD$ such that $H\circ F\circ H^{-1}$ is of the form of Eq. \ref{eqn:skew_torus}. Moreover, $H\circ F\circ H^{-1}$ is differentiable in the Y coordinate.
\end{proposition}
The proof is in  Section \ref{sect:proof-conjug}.

\textbf{Remark.} The assumption of ($E_{M}$) that there is an eigenvalue $1$ is not needed in this proposition.

{\bf Blenders revisited.} Bonatti and Diaz \cite{Blender2, Blender3} show that blenders can be expected near certain prototypical situations, specifically when there is an Anosov flow with a periodic orbit of some period $\tau$. In their remarks between their Theorem A and Corollary A they say that their proof of Theorem A shows that
(in our terminology) there is an open set of multi-chaotic diffeomorphisms near the time-$\tau$ map of that flow.  
That special orbit becomes a curve of fixed points for the time-$\tau$ map $F$, and we may expect that powers of $F$ have many invariant curves, some of which may consist entirely of periodic orbits while the rest are quasiperiodic curves. Hence the example they use is likely to have quasiperiodic curves, which is a curious coincidence with out paper.

The method of creating \HDC in their seminal paper is by perturbing this $F$ in a manner that creates blenders. 
Of course a slightly different perturbation could have made some of the curves quasiperiodic. Intuitively their initial map $F$ seems to be on the border of both their blenders and our quasiperiodic curve. While blenders are topologically robust,  quasiperiodic curves are not topologically robust. But smooth quasiperiodic curves are measure theoretically robust for almost every rotation number. 

\begin{SCfigure}
\includegraphics[width=0.5\textwidth]{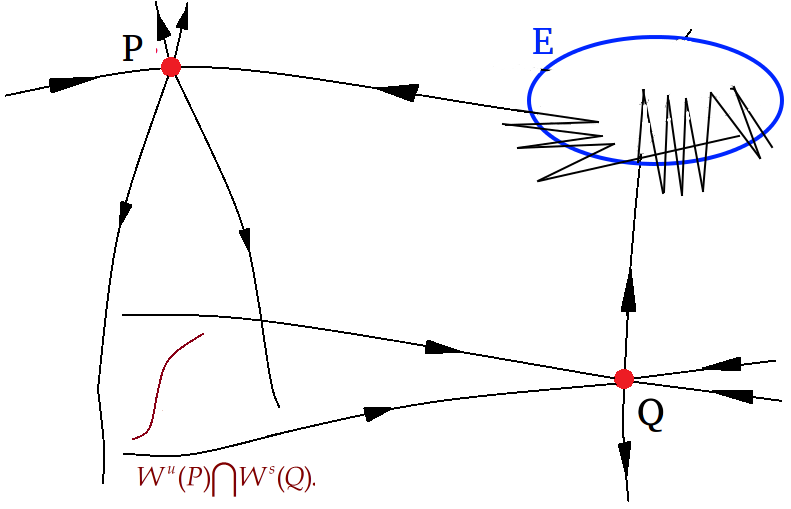}
\caption{
\textbf{Heteroclinic intersection using a quasiperiodic circle.} 
The figure provides a new simple solution to the problem of how heteroclinic cycles can be created between periodic orbits of different unstable dimensions, discussed in  Fig. \ref{fig:Unstable_dimension_illustration}. $E$ is an invariant set with a dense trajectory; and with stable and unstable dimensions both $2$, which add up to be $1$ more than the dimension of the ambient space $\mathbb{R}^3$. The stable and unstable manifolds of $P$ and $Q$ respectively are $1$-dimensional, but due to the nature of $E$ and by the $\lambda$-lemma, their closures contain the 2-dimensional manifolds $W^s(E)$ and $W^u(U)$ respectively. Thus, they have a robust intersection. 
}
\label{fig:expander}
\end{SCfigure}

The similar role that blenders and quasiperiodic curves could play has been shown in Fig. \ref{fig:expander}.  The situation is a 3-dimensional analog of the situation in  Theorem \ref{thm:Main}. In the figure, $E$ is an invariant set with a dense trajectory and with stable and unstable dimensions $2$. This forces the 1-dimensional manifolds $W^s(P)$ and $W^u(Q)$ to intersect in a neighborhood of $E$ and complete the hetero-dimensional cycle.

{\bf A connection between dense orbits and dense homoclinic points.} This paper's approach is to show that the two orbits have sets of homoclinic points that are dense in the space. The following results explains why this works and we use it in the proof of the Theorem.

\textbf{Definitions.} A map is called \textbf{topologically mixing} if for every pair of non-empty open sets $A$ and $B$, there exists some $N\in\mathbb{N}$ such that for every $n\in\mathbb{N}$ such that $n>N$, $F^{-n}(A)\cap B$ is non empty. 

\begin{proposition}\label{prop:basic_multi}
Let $\Man$ be a 2-dimensional manifold and $F:\Man\to \Man$ be a $C^1$ map. Assume there is a periodic repeller point $R$  and a periodic saddle point $\Saddle$ such that the transverse homoclinic points of each are dense in $\Man$. 
\\Then (i) the map is topologically mixing and (ii) $\Man$ is multi-chaotic.
\end{proposition}

\textbf{Proof.} The essence of the proof is that each neighborhood of a homoclinic point of a saddle or repeller contains a periodic point that is a saddle or repeller, respectively, so the hypothesis that homoclinic points are dense implies that both periodic saddle points and periodic repeller points are dense in the manifold. Recall that a point $Q$ in the unstable manifold of a repeller periodic point $R$ is a transverse homoclinic point if for some $n>0$, $F^n(Q)=R$ and $DF^n(Q)$ is nonsingular. 
F.R. Marotto introduced the term ``snapback repeller'' in 1978 \cite{SnapBR2} with improvements in \cite{SnapBR3}. 
He called a repeller periodic point $R$ a {\bf snap-back repeller} if it has a transverse homoclinic point, and showed that periodic repeller points lie in each neighborhood of each transverse homoclinic point, (and he proved much more, that snap-back repellers have ``scrambled sets'' in the ``sense of Li and Yorke'' \cite{SnapBR2}.) 

\textbf{Proof of topological mixing.} Let $A$ and $B$ be two non-empty, open subsets. Since $W^s P$ is dense in $\Man$, it intersects $A$. Since the unstable manifold $W^u$ of $P$ is dense in $\Man$, it passes through $B$. Let $z\in W^u\cap B$. 
Then there is a $\delta>0$ such that $B$ contains a $\delta$ neighborhood of $z$. 
By the lambda lemma, there is an integer $N>0$ such that for every integer $n>N$, $F^n(A)$ contains points within distance $\delta$ of $z\in W^u$ and therefore, lying within $B$. 
Therefore, for every integer $n>N$, $F^{n}(A)\cap B$ is non-empty, which implies that $F$ is topologically mixing. 

{\bf Proof that there is a dense trajectory in $\Man$.} Topological mixing implies there is a trajectory that is dense in $\Man$, so $\Man$ is multi-chaotic.

This completes the proof of  Prop. \ref{prop:basic_multi}. \qed

The proof of Theorem \ref{thm:Main} will use the following definition and lemma.

{\bf Definition.}  ${F}$ is \textbf{strongly transitive} if for every open disk $D$, there is some $N\in\mathbb{N}$ for which 
$\cup_{n=1}^{N}F^{n}(D)$ is the entire space. 

{\bf Extending to diffeomorphisms.} To keep our results as simple as possible, we have emphasized maps on $\Torus$ that are k-chaotic for $k=1$ and $2$ for the set $X=\torus$. In order to have a dense trajectory in the presence of repellers ($k=2$), the map cannot be one-to-one. The literature has emphasized diffeomorphisms, so we would like to emphasize that our methods extend to cases where there is a diffeomorphism on $\torus^3$.  Instead of the base map being $x\mapsto mx\bmod 1$ on $S^1$, it will be an Anosov diffeomorphism on $\mathbb{T}^2$, like the cat map, with an uniformly expanding direction. These systems are invertible and hyperbolic periodic points are of index 1 or 2. The unstable manifolds of index-2 saddles would intersect densely and transversely with the quasiperiodic curve, and homoclinic points of index-2 saddles would be dense. A same reasoning would apply to the stable manifolds of index-1 saddles.

\section{Proof of  Theorem \ref{thm:Main}} \label{sect:proof-Thm1}

\begin{lemma}\label{lem:strongly_trans}
Assume the hypotheses of  Theorem \ref{thm:Main}. Then ${F}$ is strongly transitive.  
In particular, 
\\(STi) if a set $X\subset\Torus$ is backward invariant (i.e. $ F^{-1}(X) \subset X$), then $X$ is dense in $\Torus$; also
\\(STii) there is a forward trajectory $(F^j(x))_{j=0}^\infty$) that is dense in $\Torus$.
\end{lemma}

\textbf{Definition.} For a map with a cone system, a differentiable curve in $\mathbb{R}^2$ or in $\Torus$ is said to be a \textbf{\textit{horizontal curve}} if its tangent at every point lies inside the expanding cone at that point on the manifold. Note that the image of a horizontal curve under the map is again a horizontal curve, with an expansion in length by a factor of at least $K$.
A \textbf{vertical curve} is a curve on the torus which is mapped under $H$ into a vertical circle of the skew-product map,
[recall that a vertical circle is represented as $S_x$, as seen in 
Eq. \ref{eqn:vertical}].

\textbf{Proof of  Lemma \ref{lem:strongly_trans}.}
Note that $F|\Gamma$ is conjugate to an irrational rotation. Let $D$ be an open disk. We will prove that $F^n(D)=\Torus$ for some integer $n>0$, thereby proving that $F$ is strongly transitive. Note that $D$ contains a horizontal curve $\lambda$. Since horizontal curves are mapped into horizontal curves by $F$ and uniformly expand by a factor of $K$, for sufficiently large $n_1\in\mathbb{N}$, $F^{n_1}(D)$ stretches across in the direction of $v_m^L$ and in particular, intersects $\Gamma$ in an open interval $I$. 
Since $F|\Gamma$ is a quasiperiodic map, for sufficiently large $n_2\in\mathbb{N}$, $\cup_{k=1}^{n_2} F^k(I)=\Gamma$.
Therefore, $\cup_{k\in[n_2]}F^{n_1+k}(I)$ contains a neighborhood $V$ of $\Gamma$. Since $\Gamma$ is compact and transverse to the cones, for sufficiently large $n_3\in\mathbb{N}$, $F^{n_3}(V)=\Torus$. 
Hence, the union of the first $n=n_1+n_2+n_3$ iterates of $D$ cover $\Torus$, proving $F$ is strongly transitive.

The proofs of (STi) and (STii) are standard, but are included for completeness. 	(STi) follows from the fact that if a backward invariant set $X$ is not dense, then the complement $U$ of the closure of $X$ is a non-empty open set that is positively invariant ($F(U) \subset U),$ which contradicts strong transitivity.
	
	(STii) follows from the following generalization, which is easy to prove. Let $(U_j)_{j=0}^\infty$ be any sequence of open sets. 
Then there is an increasing sequence of integers $(m_j)$ such that $ F^{d_j}(U_j)=\Torus$ where $d_j=m_{j+1}-m_j$. 	Writing $\bar U$ for the closure of $U$, it follows that there is a point $x_0$ for which $F^{m_j}(x_0)\subset \bar U_j$. It follows that there is a point $x_0$ so that $F^{m_j}(x_0)$ is in the closure of $U_j$ for all $j$. Choose the sequence of open sets $(U_j)$ so that every point of $\Torus$ is in infinitely many of the $U_j$ and so that the diameter of $U_j \to0$ as $j\to\infty$. Then the set $\{F^{m_j}(x_0)\}$ is dense in $\Torus$. That completes the proof of Lemma \ref{lem:strongly_trans}.

\begin{figure}
\centering
\includegraphics[width=0.7\textwidth]{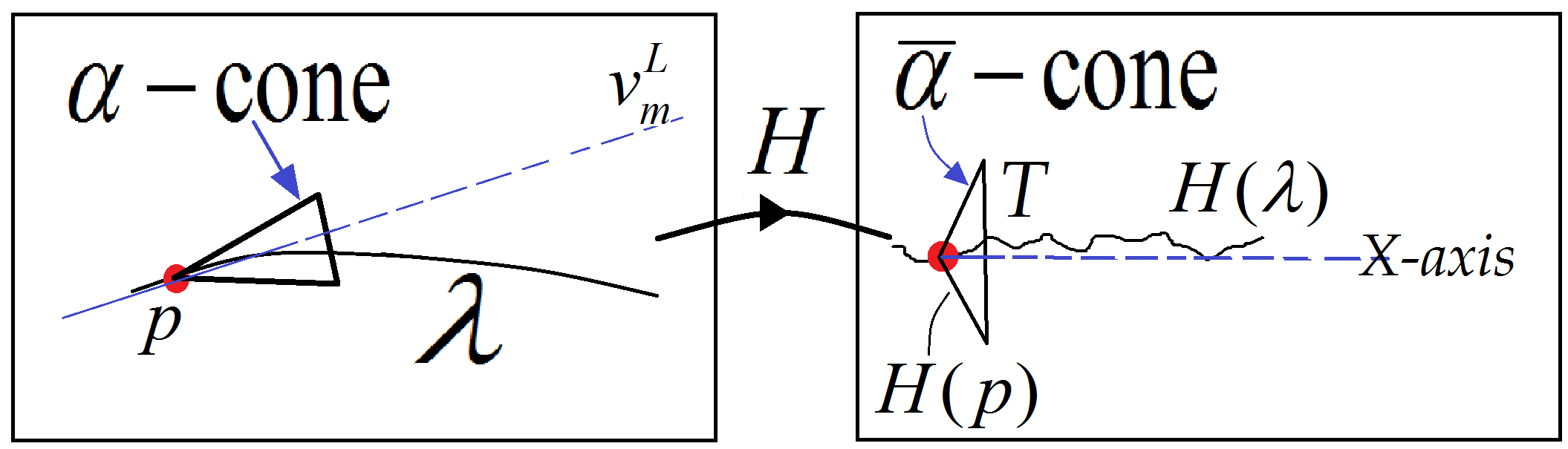}
\caption{\textbf{Illustration of Claim 1a.} The figure shows how smooth, horizontal curves on $\TorusD$ could be deformed under the conjugacy map $H$ in  Theorem \ref{thm:conj-exp-cone}. The tangent to a point $p$ on a smooth, horizontal curve $\lambda$ is always located in the $\alpha$-cone drawn at $p$. Under the mapping $H$, $\lambda$ is not necessarily $C^1$, but locally at $p$, the curve $H(\lambda)$ stays within an $\bar{\alpha}$-cone centered around the X-direction. The constant $\bar{\alpha}$ does not  depend on $\lambda$ or $p$. The triangle $T$ marked here is the triangle mentioned in Claim 1a.}
\label{fig:img_curv_H}
\end{figure}

\textbf{Some simplifying remarks.}  We are now ready to prove Theorem \ref{thm:Main}. For some power $p$, the saddle $\Saddle$ and repeller $\Repel$ are fixed points of $F^p$, and the curve $\Gamma$ is mapped onto itself. Since $F^p$ satisfies the hypotheses and since $F^p$ is multi-chaotic iff $F$ is multi-chaotic, we can assume without loss of generality that $p=1$. Therefore, $\Saddle$ and $\Repel$ are fixed points of $F$ and $\Gamma$ is mapped to itself. 

\textbf{Outline of the proof of  Theorem \ref{thm:Main}.} We will first prove that the saddle's stable and unstable manifolds are dense  in $\Torus$ in Claims 1 and 2 below. The stable and unstable manifolds of the saddle are transverse everywhere, so the transverse homoclinic points of the saddle are dense.

Next, we will prove that the repeller's stable manifold is dense in $\Torus$ and the unstable manifold is the whole $\Torus$, in Claims 3 and 4 below. Hence the homoclinic points of the repeller are also dense, and transverse as the Jacobian is everywhere non-singular. This satisfies the hypothesis of Prop. \ref{prop:basic_multi} and the conclusion of multi-chaos in $\Torus$ follows. 

\textbf{Claim 1:} $W^u(\Saddle)$ is dense in $\Torus$. 
Note that $W^u(\Saddle)$ is a horizontal curve. Differentiable horizontal curves may not be mapped under $H$ into differentiable curves, but we have the following property.

\textbf{Claim 1a:} There is a constant $\bar{\alpha}>0$ for which the following is true. Let $\lambda$ be a $C^1$ horizontal curve of finite length, and $\bar{\lambda}$ be the image $H(\lambda)$. Let $p\in \bar{\lambda}$, and $T$ be a triangle with vertex at $p$, sides at angles 
$\pm\tan^{-1}(\bar{\alpha})$ with the X axis. Then $\bar{\lambda}$ does not intersect the sides of $T$ locally at $p$. [See triangle in  Fig. \ref{fig:img_curv_H}]. (In other words, $\bar{\lambda}$ may not be differentiable but like $\lambda$, it remains bounded between ``cones''.)

\textbf{Proof.} One of the consequences of  Theorem \ref{thm:conj-exp-cone} is that $F_0$ is differentiable along the Y-direction. Therefore, the pre-images of vertical circles under $H$ will be $C^1$ curves. Let $\beta$ be the maximum expansion of $F$ along these curves. Since these curves are uniformly transverse to the invariant cones, $\beta<\infty$. Therefore, the constant $\bar{\alpha}>0$ exists, proving Claim 1a. \qed

The proof of the density of $W^u$ will be by contradiction. 
{\bf So suppose $W^u$ is not dense in $\Torus$}, i.e., there exists an open, non-empty set \boldmath $U$ \unboldmath such that $U$ and $W^u$ are disjoint. Let $U$ be the largest such open set. Since $W^u$ is forward invariant, for every integer $n$, $F^{-n}(U)\cap W^u=\emptyset$, so $U$ is backward-invariant. We will now look at the image of the components of $U$ under the conjugacy-map $H$ from  Theorem \ref{thm:conj-exp-cone}. 
By  Lemma \ref{lem:strongly_trans}, $U$ must be an open set disjoint from $W^u$ and dense in $\Torus$. Let \boldmath $K$ \unboldmath be the closure of $W^u$. Then every point on $K$ lies on the boundary of $U$.
 Note that $W^u$ is a horizontal curve and therefore, must intersect $\Gamma$ at least one point $z$. Since $F|\Gamma$ is conjugate to an irrational rotation, the orbit of $z$ must be dense in $\Gamma$. However, the orbit of $z$ is also a part of $W^u$. Therefore, $\Gamma\subset K$. 

The image of the quasiperiodic circle $\Gamma$ under $H$ is a vertical circle \boldmath $\bar\Gamma.$ \unboldmath Let \boldmath $x_\Gamma$ \unboldmath be the X coordinate of $\bar\Gamma$. Let $\lambda$ be a vertical curve contained in the open sets $U$,  such that $\bar\lambda$:=$H(\lambda)$ is at distance $\delta$ from $\bar\Gamma$.

\textbf{Claim 1b:} If $W^u$ is not dense, then $L/\delta \leq 2\bar\alpha$.

\begin{figure}
\includegraphics[width=0.32\textwidth]{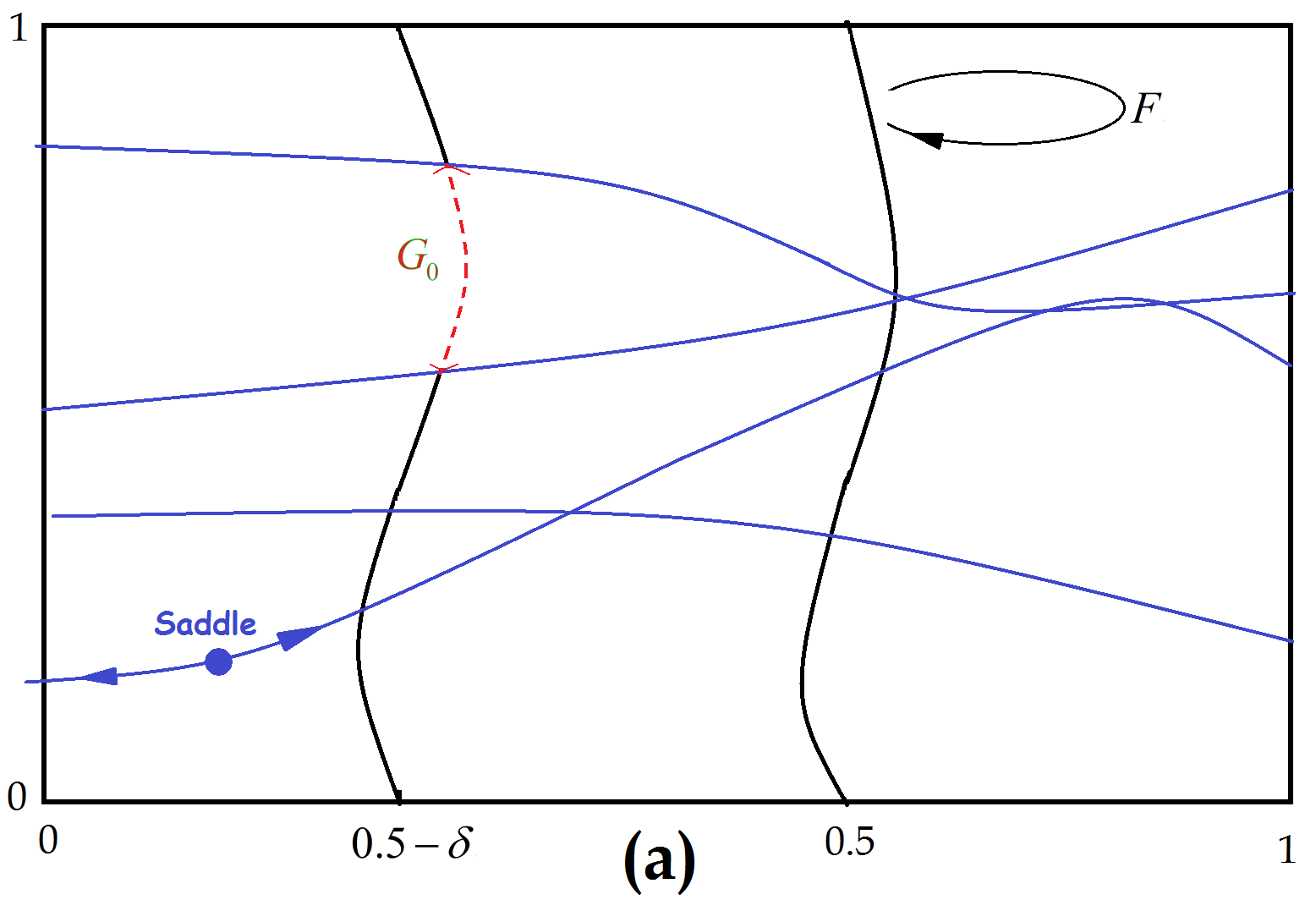}
\includegraphics[width=0.32\textwidth]{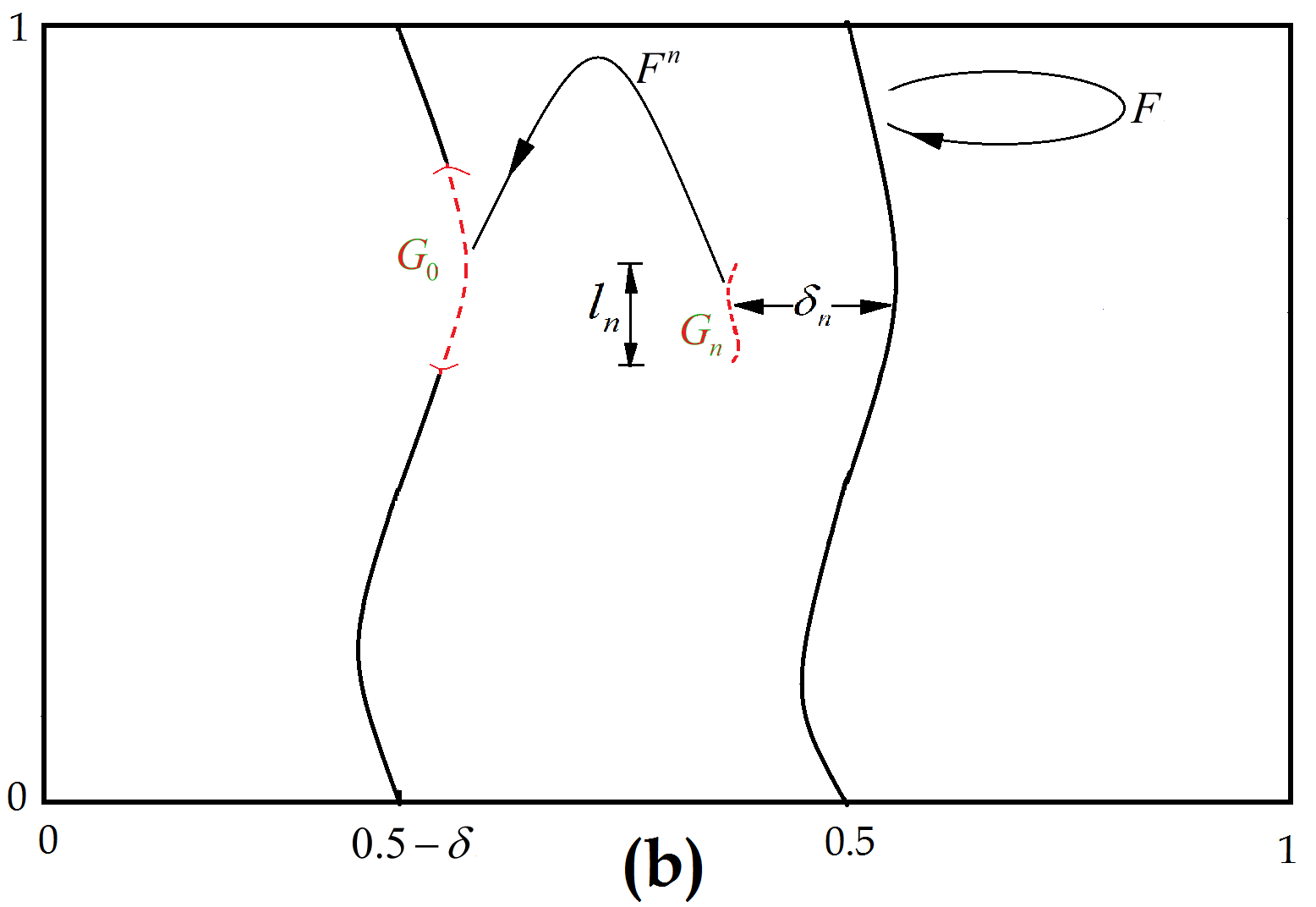}
\includegraphics[width=0.32\textwidth]{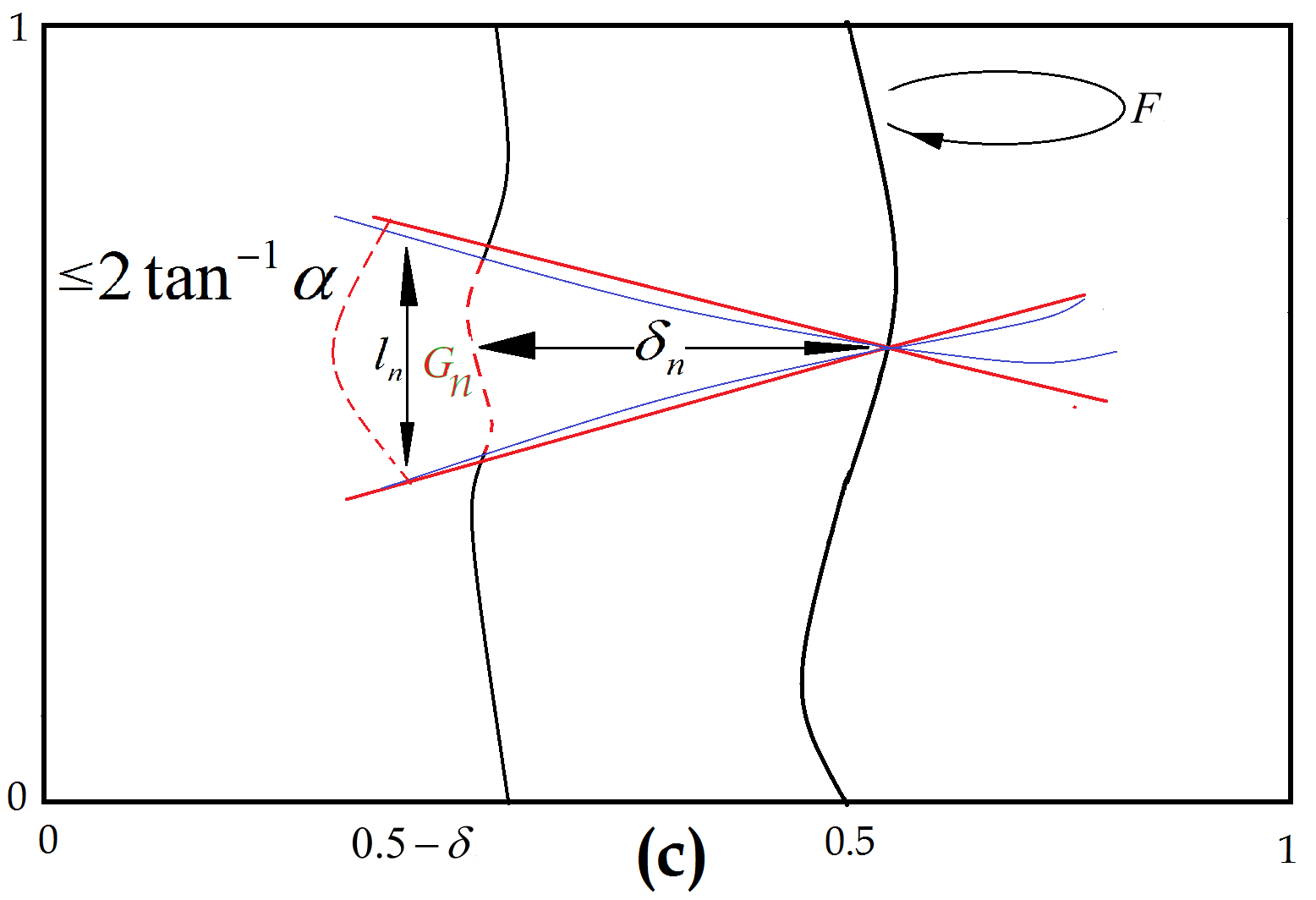}
\caption{\textbf{Density of the unstable manifold.} This illustration shows the main argument in the proof of the claim made in  Theorem \ref{thm:Main} that the unstable manifold $W^u$ is dense in $\Torus$. 
The proof is by contradiction and starts with the assumption that $W^u$ is not dense. 
The arcs marked in green in (a) are gaps, that is, open intervals in the vertical circles (shown as dashed green lines), gaps that are disjoint from $W^u$. 
The blue curves are some segments of the unstable manifold $W^u$. Fig. (b) shows an inverse image $G_{n}$ of $G_0$ under $F^n$ and it is shown in the proof that the ratio $\frac{L}{\delta}$ diverges to $\infty$ as $n\rightarrow\infty$, where $L$ and $\delta$ have been marked in the figure. In Fig. (c), we look at an arbitrary gap (green) which is part of some circle (black) with slope bounded below by some $\alpha'>\alpha$. Claim 1b gives an upper bound for the ratio $\frac{L}{\delta}$, which contradicts the fact that this ratio is unbounded.}
\label{fig:Proof_illustration}
\end{figure}

 Fig. \ref{fig:Proof_illustration} illustrates the proof. 
Let $x_1$ be the X coordinate of the curve $\bar\lambda$. Let $(x_1,y_0)$ and $(x_1,y_2)$ be the two endpoints of $\bar\lambda$, with $y_0<y_1$. Then the vertical mid-point is $(x_1,y_1)$ for $y_1=0.5(y_0+y_2)$. Then note that by assumption, (i) $|x_\Gamma-x_1|=\delta$; and (ii) the point $(x_\Gamma,y_1)$ is a point on $\bar\Gamma$. 

Since $W^u$ has dense intersections with $\bar\Gamma$, for every $\epsilon>0$, there is a point $(x_\Gamma,y_1')\in \bar\Gamma\cap W^u$ at a distance less than $\epsilon$ from $(x_\Gamma,y_1)$. Let the segment of $W^u$ starting at $(x_\Gamma,y_1')$ intersect $\bar\lambda$ at $(x_1,y_3)$. Since this curve cannot intersect with $\bar\lambda$, it intersects outside the interval $(y_0,y_2)$. Without loss of generality, let $y_3\leq y_0$. 

By Claim 1a, $|y_3-y_1|\leq|y_3-y_1'|+|y_1-y_1'| \leq \bar\alpha \delta+\epsilon$. But $L=2|y_0-y_1|\leq 2|y_3-y_1|\leq 2\alpha|\delta|+2\epsilon$. Since this inequality holds for every $\epsilon>0$, it must hold for $\epsilon=0$. This gives the inequality $\frac{L}{\delta}\leq2\bar\alpha$, proving Claim 1b.

\textbf{Claim 1c:} If $W^u$ is not dense, then $L/\delta$ is unbounded. [Compare with Claim 1b].

\textbf{Proof.}
The intersection of the open set $U$ with the vertical circle $S_{x_\Gamma+\delta}$ is a disjoint collection of maximal, open arcs. 
Let $I_0$ be such an arc with Y-span $l>0$. 
For every integer $n\geq 1$, there is a vertical arc $I_{n}$ in $S_{x_\Gamma+m^{-n}\delta}$ such that $F^n(I_{n})=I_0$. 
Let $m_y:=\max_{z\in\Torus}|\partial_y (HFH^{-1})|$. 
So $I_{n}$ lies at distance $m^{-n}\delta$ and has a Y-span of at least $\frac{l}{m_y^n}$. 
Note that $I_n$ lies in $U$, since $U$ is backward invariant.
However, the ratio of the Y-span to the distance from its mid-point is $l\left(\frac{m}{m_y}\right)^n$ which $\to\infty$ as $n\rightarrow\infty$, since $m_y<m$.   So Claim 1c is proved. See  Fig. \ref{fig:Proof_illustration} for an illustration. \qed

Therefore, if $W^u$ is assumed to be not dense, then $L/\delta$ is both bounded by Claim 1b and unbounded by Claim 1c. 
Hence $W^u$ must be dense in $\Torus$. This completes the proof of Claim 1.  \qed

\textbf{Claim 2:} $W^s(\Saddle)$ is dense in $\Torus$. Note that $W^s(\Saddle)$ will contain a portion $J$ of the vertical circle on which it lies. Then $W^s(\Saddle)$ can be constructed as the backward-invariant set $\cup_{n=0}^{\infty} F^{-n}(J)$. By  Lemma \ref{lem:strongly_trans} (STi), $W^s(\Saddle)$ is dense in $\Torus$.

\textbf{Claim 3:}  $W^s(\Repel)$ is dense in $\Torus$. Note that $W^s(\Repel)$ is the collection of all points which eventually map into $\Repel$. It is the union $\underset{n\in\mathbb{N}}{\cup}F^{-n}(\Repel)$ and is a backward-invariant set. By  Lemma \ref{lem:strongly_trans}, $W^s(\Repel)$ is dense in $\Torus$.

\textbf{Claim 4:} $W^u(\Repel)=\Torus$. Note that $\Repel$ has a repelling neighborhood $U$ and by  Lemma \ref{lem:strongly_trans}, $U$ covers $\Torus$ under a finite number of iterations of $F$. Therefore, $W^u(\Repel)=\Torus$, as claimed.

As mentioned at the beginning of the proof, establishing these claims completes the proof of  Theorem \ref{thm:Main}.

\section{Proof of  Theorem \ref{thm:conj-exp-cone}}\label{sect:proof-conjug}

\textbf{Outline of the proof.} To construct the conjugacy, 
we define a map $\upPhi:\TorusD\to S^1$ 
below which is used to define a semi-conjugacy 
$\lowPhi(z\bmod 1):=\upPhi(z)\bmod 1$, 
to the map in 
Eq. \ref{eqn:expnd_circ} below. 
Under the additional assumption of the cone system, Step 2 shows how $\lowPhi$ can be used to construct a conjugacy $H$. 
It turns out that $H$ is a continuous map but usually not differentiable. 
Step 3 shows that however, $H\circ F\circ H^{-1}$ is differentiable with respect to the Y-variable. 

{\bf Step 1: The projection semi-conjugacy \cite{Bylnd}. }
By assumption, $M$ has a left eigenvector $v_m^L$ corresponding to the eigenvalue $m$. Since $M$ is integer-valued, $v_m^L$ can be chosen to have integer entries. Note that $( v_m^L ) M^n = m^n( v_m^L )$ and $v^L_m\cdot x\in\mathbb{Z}$ for every $ x\in\mathbb{Z}^d$. Let \boldmath $k$ \unboldmath $:=\min$\{$|v_m^L x|$ : $x\in\mathbb{Z}^d$, $v_M^L x\neq 0$\}. Then note that $k$ is an integer and $k$ divides $v_M^L x$ for every $x\in\mathbb{Z}^d$.

Let $\upPhi:\mathbb{R}^d\to\mathbb{R}$ and $\lowPhi:\TorusD\to\TorusD$ be defined as follows.
\begin{equation}\label{eqn:Phi_def}
\begin{split}
\upPhi(x) &:=k^{-1}\lim_{n\to\infty}m^{-n} v_m^L\hat{F}^n(x).\\
\lowPhi(z) &=\lowPhi(z\bmod 1) :=\upPhi(z)\bmod 1.\\
\end{split}
\end{equation}

These maps were constructed by Boyland \cite{Bylnd} and defined in more generality by Franks \cite{Franks1}. The maps are well-defined, continuous and satisfy the following.
\begin{equation}\label{eqn:factoring_phi}
\begin{split}
\upPhi\circ\hat{F}(z)&=m\upPhi(z).\\
\lowPhi(F(z))&=m\lowPhi(z)\mod~1.
\end{split}
\end{equation}
\begin{equation}\label{eqn:ReverseLip_Phi}
\left||\upPhi(z_1)-\upPhi(z_2)|-|v^L_m(z_1-z_2)|\right|\leq \frac{2}{m-1}\| v_m^L \|\|G\|,\ \mbox{ for every } z_1, z_2\in\mathbb{R}^2.
\end{equation}
\begin{lemma}\label{lem:onto}[From Eq. 2.3, \cite{Bylnd}]
For every $x_0\in\mathbb{R}$ and every line \boldmath $L$ \unboldmath in $\mathbb{R}^d$ that is parallel to $ v_m^L $, there exists a point $(x,y)$ on the line for which $\upPhi(x,y)=x_0$. In particular, $\upPhi$ is onto.
\end{lemma}

For $\theta_0\in S^1$, we say the set  $\lowPhi^{-1}(\theta_0)$  is a {\bf fiber} of $\lowPhi$. Using the inequalities and  Lemma \ref{lem:onto} above, we can make the following crucial observation.

\begin{lemma}\label{prop:phi_fiber_C0}
The fibers of $\lowPhi$ are homeomorphic to $\TorusDmO$.
\end{lemma}

A map $G:N\to N$ will be called a \textbf{factor} of a map $F:M\to M$ if there is a continuous map $\pi:M\to N$ such that for every $x\in M$, $G(\pi(x)) = \pi(F(x))$.

{\bf Proof of  Prop. \ref{prop:phi_fiber_C0}.}
To prove this, we will first prove the analogous statement for $\upPhi$, i.e., for every $x_0\in\mathbb{R}$, $\upPhi^{-1}(x_0)$ is homeomorphic to the $d-1$ dimensional hyperplane $\mathbb{R}^{d-1}$. 
We will then use the fact that $\lowPhi$ is a factor of $\upPhi$ to prove the claim of this proposition. We will prove this via three claims.

\textbf{Claim A:} Two distinct points in $\upPhi^{-1}(x_0)$ cannot be connected by a horizontal curve in $\mathbb{R}^d$. To see this, assume the contrary, that there are distinct points $z_1, z_2\in\upPhi^{-1}(x_0)$ and $\gamma$ is a horizontal curve joining $z_1$ and $z_2$. 
Then for every $n\in\mathbb{N}$, $F^n(\gamma)$ is again a horizontal curve whose endpoints are $F^n(z_1)$ and $F^n(z_2)$, both lying in the fiber $\upPhi^{-1}(m^n x_0)$ by Eq. \ref{eqn:factoring_phi}. 
By Eq. \ref{eqn:ReverseLip_Phi}, we can conclude that $|v^L_m(F^n(z_1)-F^n(z_2))|\leq\frac{2}{m-1}\| v_m^L \|\|G\|$. 
Let $l$ be the length of $\gamma$. Then the length of $F^n(\gamma)$ is at least $K^n l$. 
Moreover, there is a uniform constant $\tau>0$ such that for each horizontal curve of length $l$ joining two points $A$ and $B$, 
\[|\phi(A-B)|= v_m^L \cdot(A-B)\geq\| v_m^L \|\tau l.\] 
Therefore, 
\[\frac{2}{m-1}\| v_m^L \|\|G\| \geq \phi(F^n z_1-F^n z_2)| \geq \| v_m^L \|lk^n\tau.\]
This inequality holds for every $n=1,2,\cdots$. But while the left-hand side remains bounded, the right-hand side diverges to $\infty$ as $n\rightarrow\infty$, a contradiction, proving Claim A.

\textbf{Claim B:} $\upPhi^{-1}(x_0)$ is homeomorphic to $\mathbb{R}^{d-1}$. To see this, note that every straight line parallel to $ v_m^L $ is a horizontal curve. Therefore, by Lemma \ref{lem:onto}, $\upPhi^{-1}(x_0)$ intersect every line parallel to $ v_m^L $. This combined with Claim A implies that $\upPhi^{-1}(x_0)$ intersect every line parallel to $ v_m^L $ at a unique point. Since $\upPhi$ is continuous, $\upPhi^{-1}(x_0)$ is a closed set. The family of straight lines parallel to $ v_m^L $ has a one-to-one correspondence with $\mathbb{R}^{d-1}$. Therefore, $\upPhi^{-1}(x_0)$ is homeomorphic to $\mathbb{R}^{d-1}$, proving Claim B.

\textbf{Claim C:} $\lowPhi^{-1}(\theta_0)$ is a $(d-1)$-dimensional torus. Let $x_0\in\mathbb{R}$ be some lift of $\theta_0$ under the map \emph{mod1}. Then $\lowPhi^{-1}(\theta_0)$ is the image under the map \emph{mod1} of the sets $\upPhi^{-1}(x_0+n)$, where $n$ ranges over all integers. By Claim B, each of the sets $\upPhi^{-1}(x_0+n)$ is homeomorphic to $\mathbb{R}^{d-1}$. Since $\upPhi$ is periodic and $F$ has the cone system, the images under the map \emph{mod1} of all the $(d-1)$-dimensional surfaces $\upPhi^{-1}(x_0+n)$ are the same $(d-1)$-dimensional torus, proving Claim C and thus Prop. \ref{prop:phi_fiber_C0}.

\textbf{Step 2: Constructing the conjugacy $H$.} Now that we have the factor map $\lowPhi:\TorusD\to S^1$, we can define a continuous map $H:\TorusD\to\TorusD$  which will serve as the conjugacy to a skew product. We will first introduce some notions associated with assumption ($E_{M})$.

\textbf{A tiling of the torus.} Given the integer vector $ v_m^L \in\mathbb{R}^d$, we can find another $d-1$ linearly independent integer vectors $w_1,\cdots,w_{d}$ which satisfy the following.
\\(i) For each $i=1\cdots,d-1$, $w_i$ is perpendicular to $ v_m^L $.
\\(ii) The volume of the parallelepiped \boldmath $\Tile$ \unboldmath formed by \{$w_1,\cdots,w_{d}$\} have unit volume.
\\(iii) $\phi(z)=0$ on the side of $\Tile$ which contains $w_1,\cdots,w_{d-1}$, and is $1$ on the opposite face of $\Tile$.

To see this, first let $\Plane$ be the $(d-1)$-hyperplane of vectors perpendicular to $ v_m^L $. So we can pick $d-1$ linearly independent integer vectors $w_1,\cdots,w_{d-1}$ in $\Plane$ such that the $(d-1)$-dimensional parallelepiped whose edges include 
 $w_1,\cdots,w_{d-1}$ has no lattice point in its interior or in the interior of any of its faces. There are two lattice points away from $\Plane$ which are closest to $\Tile$. Let $w_d$ be the one among them for which $\phi(w_d)>0$. Let $\Tile$ be the parallelepiped 
whose edges include $v_m^L,w_1,\cdots,w_{d}$.

Now note that the $\Tile$ has no lattice point in its interior or in the interior of any of its faces. Therefore, by the $d$-dimensional Pick's formula, $\Tile$ has volume 1. Therefore, $\Tile$ forms a periodic tiling of $\mathbb{R}^d$ and $\Tile\mod Z$  is homeomorphic to  $\TorusD$.

We claim that $\phi(w_d)=1$. To see this, first note that by construction, the two planes $\Plane$ and $w_d\oplus\Plane$ have no lattice points inside them and $\phi$ is constant on each translate of $\Plane$. Secondly, $\phi$ is integer valued on every translate of $\Plane$ which passes through a lattice point and hence is integer valued on $w_d\oplus\Plane$. In particular, $\phi(w_d)\geq 1$. Moreover, since $ v_m^L $ has no common denominator for its entries, there is a lattice point $z_*\in\mathbb{R}^d$ for which $\phi(z_*)=1$. therefore, The plane $\phi\equiv1$, which is the same as $z_*\oplus\Plane$, cannot pass through the interior of $\Tile$. Therefore, this plane must be the ``upper'' boundary of $\Tile$ and $\phi(w_d)=1$.

\textbf{The conjugacy.} Let \boldmath $\mbox{proj}_W$ \unboldmath denote the orthogonal projection of vectors in $\mathbb{R}^d$ onto the plane spanned by $w_1,\cdots,w_d$. Let $H:\TorusD\to\TorusD$ be defined as $H(z)=(\lowPhi(z),\mbox{proj}_W z \bmod 1)$. The rest of this section will be used to prove that $H$ is the claimed conjugacy of  Theorem \ref{thm:conj-exp-cone}.

\textbf{$H$ is a homeomorphism.} Since $\lowPhi$ is continuous, $H$ is continuous. Since $\Torus$ is a compact set and $H$ is continuous, to prove that $H$ is a homeomorphism, it is enough to show that the map is both one-to-one and onto. Let $(x_0,y_0)\in S^1\times\mathbb{T}^{d-1}$. Then $\lowPhi^{-1}(x_0)$ is topologically $\mathbb{T}^{d-1}$  by  Prop. \ref{prop:phi_fiber_C0} and uniformly transverse to all lines by
Lemma \ref{lem:onto}. On the other hand, the set of points \{$z\in \TorusD$ : $proj_W(z)=y_0$\} is a circle parallel to $v_M^L$. Hence they intersect at some point $z_0$ and $H(z_0)=(x_0,y_0)$. therefore, $H$ is onto.

It also follows from  Prop. \ref{prop:phi_fiber_C0} and  Lemma \ref{lem:onto} that $H$ is also one-to-one. 

\textbf{$H$ satisfies the conjugacy relation.} For each $(x_0,y_0)\in\TorusD$, let $(x_1,y_1):=H^{-1}(x_0,y_0)$, $(x_2,y_2):=F(x_1,y_1)$ and $(x_3,y_3):=H(x_2,y_2)$. To show that $H$ is the desired conjugacy have to show that $x_3=mx_0\ (\mod~1)$. Note that $x_3=\lowPhi(x_2,y_2)=\lowPhi\circ F(x_1,y_1)$. By 
Eq. \ref{eqn:factoring_phi}, $x_3=m\times \lowPhi(x_1,y_1)$. But since $(x_1,y_1)=H^{-1}(x_0,y_0)$, $\lowPhi(x_1,y_1)$ must be equal to $x_0$. Therefore, $x_3=mx_0\ (\mod~1)$, proving that $H$ satisfies the conjugacy relation. \qed

\textbf{Step 3: Proof of smoothness in $Y$ variable.} The last claim of the Proposition, that $H\circ F\circ H^{-1}$ is $C^1$ with respect to the Y-variable, will now be proved. The claim is equivalent to saying that $H\circ F\circ H^{-1}$ is $C^1$ along the vertical circles (see Eq. \ref{eqn:vertical}). Note that the pre-images of the vertical circles are the fibers $\Phi^{-1}(\theta)$, so it is equivalent to proving that $F$ is $C^1$ along these curves. Since $F$ is $C^1$ on the whole manifold, it is enough to prove that these fibers are $C^1$ curves. This follows from Pugh, Shub and Wilkinson's technique of proving $C^1$-smoothness of central foliations in partially hyperbolic systems, see Theorem B, \cite{Holder1}. Although we do not work with a diffeomorphism or a partially hyperbolic splitting, their technique applies in our case, because the expansion along the fibers is less than the expansion along the cones. This completes the proof of  Theorem \ref{thm:conj-exp-cone}.



\section*{Acknowledgments}  
The authors thank M. Sanjuan for his extensive comments. 

\bibliographystyle{siamplain}
\bibliography{Multichaos_bibliography}
\end{document}